\documentstyle[12pt]{article}
\begin{document}
\newcommand{\text}[1]{\mbox{{\rm #1}}}
\newcommand{\gd}{\delta}
\newcommand{\itms}[1]{\item[[#1]]}
\newcommand{\nin}{\in\!\!\!\!\!/}
\newcommand{\sub}{\subset}
\newcommand{\cntd}{\subseteq}
\newcommand{\go}{\omega}
\newcommand{\Pa}{P_{a^\nu,1}(U)}
\newcommand{\fx}{f(x)}
\newcommand{\fy}{f(y)}
\newcommand{\gD}{\Delta}
\newcommand{\gl}{\lambda}
\newcommand{\gL}{\Lambda}
\newcommand{\half}{\frac{1}{2}}
\newcommand{\sto}[1]{#1^{(1)}}
\newcommand{\stt}[1]{#1^{(2)}}
\newcommand{\Z}{\hbox{\sf Z\kern-0.720em\hbox{ Z}}}
\newcommand{\singcolb}[2]{\left(\begin{array}{c}#1\\#2
\end{array}\right)}
\newcommand{\ga}{\alpha}
\newcommand{\gb}{\beta}
\newcommand{\gga}{\gamma}
\newcommand{\ul}{\underline}
\newcommand{\ol}{\overline}
\newcommand{\qed}{\kern 5pt\vrule height8pt width6.5pt depth2pt}
\newcommand{\Lrraro}{\Longrightarrow}
\newcommand{\Nb}{|\!\!/}
\newcommand{\NN}{{\rm I\!N}}
\newcommand{\bsl}{\backslash}
\newcommand{\gt}{\theta}
\newcommand{\op}{\oplus}
\newcommand{\C}{{\bf C}}
\newcommand{\Op}{\bigoplus}
\newcommand{\CR}{{\cal R}}
\newcommand{\tr}{\bigtriangleup}
\newcommand{\grr}{\omega_1}
\newcommand{\ben}{\begin{enumerate}}
\newcommand{\een}{\end{enumerate}}
\newcommand{\ndiv}{\not\mid}
\newcommand{\bab}{\bowtie}
\newcommand{\hal}{\leftharpoonup}
\newcommand{\har}{\rightharpoonup}
\newcommand{\ot}{\otimes}
\newcommand{\OT}{\bigotimes}
\newcommand{\bwe}{\bigwedge}
\newcommand{\gep}{\varepsilon}
\newcommand{\gs}{\sigma}
\newcommand{\rbraces}[1]{\left( #1 \right)}
\newcommand{\bbox}{$\;\;\rule{2mm}{2mm}$}
\newcommand{\sbraces}[1]{\left[ #1 \right]}
\newcommand{\bbraces}[1]{\left\{ #1 \right\}}
\newcommand{\OO}{_{(1)}}
\newcommand{\TT}{_{(2)}}
\newcommand{\FF}{_{(3)}}
\newcommand{\minus}{^{-1}}
\newcommand{\CV}{\cal V}
\newcommand{\CVs}{\cal{V}_s}
\newcommand{\un}{U_q(sl_n)'}
\newcommand{\on}{O_q(SL_n)'}
\newcommand{\slq}{U_q(sl_2)}
\newcommand{\olq}{O_q(SL_2)}
\newcommand{\UU}{U_{(N,\nu,\go)}}
\newcommand{\HH}{H_{n,q,N,\nu}}
\newcommand{\ct}{\centerline}
\newcommand{\bs}{\bigskip}
\newcommand{\qua}{\rm quasitriangular}
\newcommand{\ms}{\medskip}
\newcommand{\noin}{\noindent}
\newcommand{\mat}[1]{$\;{#1}\;$}
\newcommand{\raro}{\rightarrow}
\newcommand{\map}[3]{{#1}\::\:{#2}\raro{#3}}
\newcommand{\alg}{{\rm Alg}}
\def\newtheorems{\newtheorem{theorem}{Theorem}[subsection]
                 \newtheorem{cor}[theorem]{Corollary}
                 \newtheorem{prop}[theorem]{Proposition}
                 \newtheorem{lemma}[theorem]{Lemma}
                 \newtheorem{defn}[theorem]{Definition}
                 \newtheorem{Theorem}{Theorem}[section]
                 \newtheorem{Corollary}[Theorem]{Corollary}
                 \newtheorem{Proposition}[Theorem]{Proposition}
                 \newtheorem{Lemma}[Theorem]{Lemma}
                 \newtheorem{Defn}[Theorem]{Definition}
                 \newtheorem{Example}[Theorem]{Example}
                 \newtheorem{Remark}[Theorem]{Remark}
                 \newtheorem{claim}[theorem]{Claim}
                 \newtheorem{sublemma}[theorem]{Sublemma}
                 \newtheorem{example}[theorem]{Example}
                 \newtheorem{remark}[theorem]{Remark}
                 \newtheorem{question}[theorem]{Question}
                 \newtheorem{conjecture}{Conjecture}[subsection]}
\newtheorems
\newcommand{\proof}{\par\noindent{\bf Proof:}\quad}
\newcommand{\dmatr}[2]{\left(\begin{array}{c}{#1}\\
                            {#2}\end{array}\right)}
\newcommand{\doubcolb}[4]{\left(\begin{array}{cc}#1&#2\\
#3&#4\end{array}\right)}
\newcommand{\qmatrl}[4]{\left(\begin{array}{ll}{#1}&{#2}\\
                            {#3}&{#4}\end{array}\right)}
\newcommand{\qmatrc}[4]{\left(\begin{array}{cc}{#1}&{#2}\\
                            {#3}&{#4}\end{array}\right)}
\newcommand{\qmatrr}[4]{\left(\begin{array}{rr}{#1}&{#2}\\
                            {#3}&{#4}\end{array}\right)}
\newcommand{\smatr}[2]{\left(\begin{array}{c}{#1}\\
                            \vdots\\{#2}\end{array}\right)}

\newcommand{\ddet}[2]{\left[\begin{array}{c}{#1}\\
                           {#2}\end{array}\right]}
\newcommand{\qdetl}[4]{\left[\begin{array}{ll}{#1}&{#2}\\
                           {#3}&{#4}\end{array}\right]}
\newcommand{\qdetc}[4]{\left[\begin{array}{cc}{#1}&{#2}\\
                           {#3}&{#4}\end{array}\right]}
\newcommand{\qdetr}[4]{\left[\begin{array}{rr}{#1}&{#2}\\
                           {#3}&{#4}\end{array}\right]}

\newcommand{\qbracl}[4]{\left\{\begin{array}{ll}{#1}&{#2}\\
                           {#3}&{#4}\end{array}\right.}
\newcommand{\qbracr}[4]{\left.\begin{array}{ll}{#1}&{#2}\\
                           {#3}&{#4}\end{array}\right\}}

\title{The Representation Theory of Co-triangular Semisimple Hopf Algebras}
\author{Pavel Etingof
\\Department of Mathematics\\
Harvard University\\Cambridge, MA 02138
\and Shlomo Gelaki\\Department of Mathematics\\
University of Southern California\\Los Angeles, CA 90089}
\date{December 18, 1998}
\maketitle

\section{Introduction}
In [EG1, Theorem 2.1] we prove that {\em any} semisimple triangular Hopf
algebra $A$ over
an
algebraically closed field of characteristic $0$ (say the field of
complex numbers $\C$) is obtained
from a finite group after twisting the ordinary comultiplication of its
group algebra in the sense
of Drinfeld [D]; that is $A=\C[G]^J$ for some finite group $G$ and a
twist $J\in \C[G]\ot \C[G].$ In
[EG2] we show how to construct twists for
certain solvable non-abelian groups by iterating twists of their abelian
subgroups, and thus
obtain new non-trivial semisimple triangular Hopf algebras. We also
show how any non-abelian finite group which admits a bijective 1-cocycle
with
coefficients in an abelian group, gives rise to a
non-trivial semisimple minimal triangular Hopf algebra. Such non-abelian
groups (which
are necessarily solvable [ESS]) exist in abundance and were constructed
in [ESS] in connection with
set-theoretical solutions to the quantum Yang-Baxter equation.

If $A$ is minimal triangular then $A$ and $A^{*op}$ are isomorphic as
Hopf algebras. But
any non-trivial semisimple triangular $A$ which is not minimal, gives
rise to a new
Hopf algebra $A^*,$ which is also semisimple by [LR]. These are very
interesting semisimple Hopf
algebras which arise from finite groups, and they are abundant by the
constructions given in
[EG2]. Generally, the dual Hopf algebra of a triangular Hopf algebra is
called {\em co-triangular}
in the literature.

In this paper we explicitly describe the representation theory of
co-triangular
semisimple Hopf algebras $A^*=(\C[G]^J)^*$ in terms of
representations of some
associated groups. As a corollary we prove that Kaplansky's 6th
conjecture [K] holds for $A^*;$
that is that the dimension of any irreducible representation of $A^*$
divides the dimension of $A.$

We note that we have used in an essential way the results of the paper
[Mo],
from which we learned a great deal.
\section{Preliminaries}
\subsection{Projective Representations and Central Extensions}
Here we recall some basic facts about projective representations
and central extensions. They can be found in textbooks, e.g.
[CR, Section 11E].

A projective representation of a group $\Gamma$ is a vector space $V$
together with a
homomorphism of groups $\pi_{_V}:\Gamma\raro PGL(V),$ where $PGL(V)\cong
GL(V)/\C$ is the projective
linear group.

A linearization of a projective representation $V$ of $\Gamma$
is a central extension $\hat\Gamma$ of $\Gamma$ by a central subgroup
$\zeta$ together with a linear representation
$\tilde\pi_V:\hat\Gamma\to GL(V)$ which descends to $\pi_{_V}$.
If $V$ is a finite-dimensional projective representation of $\Gamma$
then there exists a linearization of $V$ such that $\zeta$ is finite
(in fact, one can make $\zeta=\Z/(\text{dim}V)\Z$).

Any projective representation $V$ of $\Gamma$ canonically defines a
cohomology class
$[V]\in H^2(\Gamma,\C^*)$. The representation $V$ can be lifted to
a linear representation of $\Gamma$ if and only if $[V]=0$.
\subsection{The Algebras Associated With a Twist}
Let $H$ be a finite group, and let $J\in \C[H]\ot \C[H]$ be a {\em
minimal} twist
(see [EG1]). That is, the right (and left) components of the
R-matrix $R=J_{21}^{-1}J$ span $\C[H].$ Define two coalgebras
$(A_1,\Delta_1,\varepsilon),(A_2,\Delta_2,\varepsilon)$ as
follows: $A_1=A_2=\C[H]$ as vector spaces, the coproducts are determined
by
$$\Delta_1(x)=(x\ot x)J,\;\Delta_2(x)=J^{-1}(x\ot x)$$
for all $x\in H,$ and $\varepsilon$ is the ordinary counit of $\C[H].$
Note that
since $J$ is a twist, $\Delta_1$ and $\Delta_2$ are indeed
coassociative. Clearly the dual
algebras $A_1^*$ and $A_2^*$
are spanned by $\{\delta_h|h\in H\},$ where $\delta_h(h')=\delta_{hh'}.$
\begin{theorem}\label{movshev}
Let $A_1^*$ and $A_2^*$ be as above. The following hold:
\ben
\item $A_1^*$ and $A_2^*$ are $H-$algebras via
$$\rho_1(h)\delta_y=\delta_{hy},\;\;\rho_2(h)\delta_y=\delta_{yh^{-1}}$$
respectively.
\item $A_1^*\cong A_2^{*op}$ as $H-$algebras (where $H$ acts on
$A_2^{*op}$ as it does on $A_2^*$).
\item
The algebras $A_1^*$ and $A_2^*$ are simple, and are isomorphic as
$H-$modules to the regular
representation $R_H$ of $H.$
\een
\end{theorem}
\proof The proof of part 1 is straightforward.

The proof of part 3 follows from the results in [Mo]. Namely, it follows
from
[Mo, Proposition 14] that in the case of a minimal twist the group $St$
defined in [Mo] (which is,
by
definition, a subgroup of $H$) coincides with $H$. Therefore,
by [Mo, Propositions 6,7] the algebras $A_1^*$ and $A_2^*$ are simple.
Furthermore, by [Mo, Propositions 11,12], the actions of $H$ on $A_1^*$
and $A_2^*$ are
isomorphic to the regular representation of $H$.

Let us prove part 2. Let
$S_0,\Delta_0,m_0$ denote the standard antipode, coproduct and
multiplication of $\C[H],$ and
define $Q=m_0(S_0\ot I)(J).$ Then it is straightforward to
verify that $Q$ is invertible, and $(S_0\ot S_0)(J)=(Q\ot
Q)J_{21}^{-1}\Delta_0(Q)^{-1}$ (see e.g. (2.17) in [Ma, Section 2.3]).
Hence the map $A_2^*\raro A_1^{*op},$
$\delta_x\mapsto \delta_{S_0(x)Q^{-1}}$ determines an $H-$algebra
isomorphism. \qed
\begin{cor}\label{square}
Let $A$ be a semisimple minimal triangular Hopf algebra over $\C$ with
Drinfeld element $u.$ If
$u=1$ then $dimA$ is a square, and if $u\ne 1$ then $2dimA$ or $dim A$
is a square.
\end{cor}
\proof The first statement follows from [EG1, Theorem 2.1] and part 3 of
Theorem \ref{movshev}. To
prove the second statement, let $(A,R,u)$ be a semisimple minimal
triangular Hopf algebra with
$u\ne 1$, and $(A,R',u')$ be obtained from $(A,R,u)$ by changing $R$ so
that the new Drinfeld
element $u'=1$ (as in [EG1]). Then $(A,R')=({\bf C}[H]^J,J_{21}^{-1}J)$
for some finite group $H.$
Let $A_{min}=\C[H']^J$ be the minimal triangular Hopf subalgebra of
$(A,R')$ where $H'\subset H$ is
a subgroup, and $J$ is a minimal twist for $H'.$ It is clear that $H$ is
generated by $H'$ and
$u.$ Since $u$ is a central grouplike element of order $2,$ we get that
the index of $H'$ is at
most $2.$ This implies our statement.\qed

Let $A_1^*,A_2^*$ be the $H$-algebras as in Theorem \ref{movshev}.
Since the algebras $A_1^*$, $A_2^*$ are simple, the actions of $H$ on
$A_1^*,$ $A_2^*$
give rise to projective representations $H\raro PGL(|H|^{1/2},\C).$
We will denote these projective representations by $V_1,$
$V_2$ (they can be thought of as the simple modules over $A_1^*,$
$A_2^*$, with the induced
projective action of $H$). Note that part 2 of Theorem \ref{movshev}
implies that $V_1,$ $V_2$ are
dual to each other, hence that $[V_1]=-[V_2].$
\section{The Main Result} Let $(A,R)$ be a semisimple triangular Hopf
algebra over $\C,$ and assume
that the Drinfeld
element $u$ is $1$
(this can be always achieved by a simple modification of $R$,
without changing the Hopf algebra structure [EG1]). Then
by [EG1, Theorem 2.1], there exist finite groups
$H\subset G$ and a minimal twist $J\in \C[H]\otimes \C[H]$ such that
$(A,R)\cong (\C[G]^J,J_{21}^{-1}J)$ as triangular Hopf algebras.
So from now on we will assume that $A$ is of this form.

Consider the dual Hopf algebra $A^*$. It has a basis of
$\delta$-functions $\delta_g$.
The first simple but important fact about the structure of $A^*$ as an
algebra is:
\begin{Proposition}\label{propo1}
Let $Z$ be a double coset of $H$ in $G,$ and
$\displaystyle{A_Z^*=\oplus_{g\in
Z}\C\delta_g\subset A^*}.$ Then $A_Z^*$ is a subalgebra of $A^*,$ and
$\displaystyle{A^*=\oplus_Z
A_Z^*}$ as algebras.
\end{Proposition}
\proof Straightforward. \qed

Thus, to study the representation theory of $A^*,$ it is sufficient to
describe the representations
of $A_Z^*$ for any $Z.$

Let $Z$ be a double coset of $H$ in $G$, and let $g\in Z$.
Let $K_g=H\cap gHg^{-1},$ and define the embeddings
$\theta_1,\theta_2:K_g\to H$
given by $\theta_1(a)=g^{-1}ag$, $\theta_2(a)=a$. Denote by $W_i$ the
pullback
of the projective $H$-representation $V_i$ to $K_g$ by means of
$\theta_i$, $i=1,2$.

Our main result is the following theorem, which is proved in the next
section.
\begin{Theorem}\label{main}
Let $W_1,W_2$ be as above, and let $(\hat K_g,\tilde\pi_{_W})$ be any
linearization of the projective representation
$W=W_1\otimes W_2$ of $K_g.$ Let $\zeta$ be the kernel of the projection
$\hat K_g\to K_g,$ and $\chi:\zeta\to \C^*$ be the character by which
$\zeta$ acts in $W$. Then there exists a 1-1 correspondence between
isomorphism classes of
irreducible representations of $A_Z^*$ and isomorphism classes of
irreducible representations of
$\hat K_g$ with $\zeta$ acting by $\chi$.
If a representation $Y$ of $A_Z^*$ corresponds to a representation
$X$ of $\hat K_g$ then $\displaystyle{dimY=\frac{|H|}{|K_g|}dimX}.$
\end{Theorem}

As a corollary we get Kaplansky's 6th conjecture [K] for semisimple
co-triangular Hopf algebras.
\begin{Corollary}\label{kap}
The dimension of any irreducible representation of a semisimple
co-triangular Hopf algebra divides
the dimension of the Hopf algebra.
\end{Corollary}
\proof Since $dimX$ divides $|K_g|$ (see e.g. [CR, Proposition 11.44]),
we have that \linebreak
$\displaystyle{\frac{|G|}{\frac{|H|}{|K_g|}dimX}=
\frac{|G|}{|H|}\frac{|K_g|}{dimX}}$ and the result
follows. \qed

In some cases the classification of representations of $A_Z^*$ is even
simpler.
Namely, let $\ol{g}\in Aut(K_g)$ be given by $a\mapsto g^{-1}ag.$ Then
we have:
\begin{Corollary}\label{corf}
If the cohomology class
$[W_1]$ is $\ol{g}-$invariant then irreducible representations of
$A_Z^*$ correspond in a 1-1
manner to irreducible representations of
$K_g$, and if $Y$ corresponds to $X$ then
$\displaystyle{dimY=\frac{|H|}{|K_g|}dimX}.$
\end{Corollary}
\proof
For any $\ga \in Aut(K_g)$ and $f\in Hom((K_g)^n,\C^*),$ let $\ga \circ
f\in Hom((K_g)^n,\C^*)$
be given by $(\ga \circ f)(h_1,\dots,h_n)=f(\ga (h_1),\dots,\ga (h_n))$
(which determines the
action of $\ga$ on $H^i(K_g,\C^*)).$ Then it follows from the identity
$[V_1]=-[V_2],$ given at
the end of Section 2,
that $[W_1]=-\ol{g}\circ[W_2].$ Thus, in our situation $[W]=0$,
hence $W$ comes from a linear representation of
$K_g$. Thus, we can set $\hat K_g=K_g$ in the theorem, and the result
follows.\qed
\begin{Example} {\rm Let $p>2$ be a prime number, and $H=(\Z/p\Z)^2$
with the standard symplectic
form $(,):H\times H\to \C^*$ given by
$\displaystyle{\left((x,y),(x',y')\right)=
e^{2\pi i(xy'-yx')/p}.}$ Then the element
$\displaystyle{J=p^{-2}\sum_{a,b\in H}(a,b)a\ot b}$ is a minimal twist
for $\C[H].$ Let $g\in GL_2(\Z/p\Z)$ be an automorphism of $H$,
and $G_0$ be the cyclic group generated by $g$. Let $G$ be the semidirect
product of $G_0$ and $H$. It is easy to see that in this case, the
double cosets are ordinary
cosets $g^kH$, and $K_{g^k}=H$. Moreover, one can show either explicitly
or using [Mo,
Proposition 9], that $[W_1]$ is a generator
of $H^2(H,\C^*)$ which is isomorphic to $\Z/p\Z.$ The element
$g^k$ acts on $[W_1]$ by multiplication by $det(g^k)$. Therefore,
by Corollary \ref{corf},
the algebra $A_{g^kH}^*$ has $p^2$ 1-dimensional representations
(corresponding to linear representations of $H$) if $det(g^k)=1$.

However, if $det(g^k)\ne 1$, then $[W]$ generates $H^2(H,\C^*)$. Thus,
$W$ comes from a linear
representation of the Heisenberg
group $\hat H$ (a central extension of $H$ by $\Z/p\Z$) with some
central character $\chi$. Thus, $A_{g^kH}^*$ has one p-dimensional
irreducible representation, corresponding
to the unique irreducible representation of $\hat H$
with central character $\chi$ (which is $W$).}
\end{Example}
\section{Proof of Theorem \ref{main}}
Let $Z\subset G$ be a double coset of $H$ in $G,$ and let $A_1,A_2$ be
as in Subsection 2.2.
For any $g\in Z$ define the linear map
$$F_g:A_Z^*\raro A_2^*\ot A_1^*,\;\;\delta_y\mapsto \sum_{h,h'\in
H:y=hgh'}\delta_h\ot
\delta_{h'}.$$
\begin{Proposition}\label{prop4}
Let $\rho_1,\rho_2$ be as in Theorem \ref{movshev}. Then:
\ben
\item
The map $F_g$ is an injective homomorphism of algebras.
\item
$F_{aga'}(\varphi)=(\rho_2(a)\ot \rho_1(a')^{-1})F_g(\varphi)$ for any
$a,a'\in H,$ $\varphi\in
A_Z^*.$
\een
\end{Proposition}
\proof
1. It is straightforward to verify that the map $F_g^*:A_2\ot
A_1\raro A_Z$ is determined by $h \ot h'\mapsto hgh',$ and that it is a
surjective homomorphism of
coalgebras. Hence the result follows.

\noindent
2. Straightforward. \qed

For any $a\in K_g$ define $\rho(a)\in Aut(A_2^*\ot A_1^*)$ by
$\rho(a)=\rho_2(a)\ot \rho_1(a^g),$ where $a^g=g^{-1}ag$ and
$\rho_1,\rho_2$ are as in Theorem
\ref{movshev}. Then $\rho$ is an action of $K_g$ on $A_2^*\ot A_1^*.$
\begin{Proposition}\label{prop5}
Let $U_g=(A_2^*\otimes A_1^*)^{\rho(K_g)}$ be the algebra of invariants.
Then $Im(F_g)=U_g,$
so $A_Z^*\cong U_g$ as algebras.
\end{Proposition}
\proof
It follows from Proposition \ref{prop4} that
$Im(F_g)\subseteq U_g,$
and $\displaystyle{rk(F_g)=dimA_Z^*=\frac{|H|^2}{|K_g|}}.$ On the
other hand, by Theorem \ref{movshev},
$A_1^*,A_2^*$ are isomorphic to the regular representation $R_H$ of $H.$
Thus, $A_1^*,A_2^*$ are
isomorphic to $\displaystyle{\frac{|H|}{|K_g|}R_{K_g}}$ as
representations of $K_g$, via
$\rho_1(a),\rho_2(a^g)$. Thus,
$\displaystyle{A_2^*\ot A_1^*\cong \frac{|H|^2}{|K_g|^2}(R_{K_g}\ot
R_{K_g})\cong \frac{|H|^2}{|K_g|}R_{K_g}}.$ So $U_g$ has dimension
$|H|^2/|K_g|,$ and the result
follows. \qed

Now we are in a position
to prove Theorem \ref{main}.  Since $W_i\ot W_i^*\cong A_i^*$ for
$i=1,2,$ it follows from Theorem \ref{movshev} that
$\displaystyle{W_1\ot W_2\otimes W_1^*\otimes W_2^*\cong
\frac{|H|^2}{|K_g|}R_{K_g}}$ as
$\hat{K_g}$ modules. Thus, if
$\chi_{_W}$ is the character of $W=W_1\otimes W_2$ as a $\hat{K_g}$
module then
$$|\chi_{_W}(x)|^2=0,\,x\notin \zeta\;\;
and\;\;|\chi_{_W}(x)|^2=|H|^2,\,x\in\zeta.$$
Therefore,
$$\chi_{_W}(x)=0,\,x\notin \zeta\;\; and\;\;\chi_{_W}(x)=|H|\cdot
x_{_W},\,x\in\zeta,$$
where $x_{_W}$ is the root of unity by which $x$ acts in $W$.
Now, it is clear from the definition of $U_g$ (see Proposition
\ref{prop5}) that
$U_g=End_{\hat{K_g}}(W).$
Thus if $\displaystyle{W=\bigoplus_{M\in Irr(\hat{K_g})} W(M)\ot M},$
where
$W(M)=Hom_{\hat{K_g}}(M,W)$ is the multiplicity space, then
$\displaystyle{U_g=\bigoplus_{M:W(M)\ne 0}End_{\C}(W(M))}.$ So
$\{W(M)|W(M)\ne 0\}$ are the
irreducible representations of $U_g.$ Thus the following implies the
theorem:

\noindent
{\bf Lemma.}
\ben
\item
$W(M)\ne 0$ if and only if for all $x\in \zeta,$ $x_{_{|M}}=x_{_{|W}}.$
\item
If $W(M)\ne 0$ then $\displaystyle{dimW(M)=\frac{|H|}{|K_g|}dimM}.$
\een
{\bf Proof of the Lemma.} The "only if" part of 1 is clear. For the "if"
part compute $dimW(M)$ as
the
inner product $(\chi_{_{W}},\chi_{_M}).$ We have
$$(\chi_{_{W}},\chi_{_M})=\sum_{x\in
\zeta}\frac{|H|}{|\hat{K_g}|}x_{_{|W}}\cdot
dimM\cdot \bar x_{_{|M}}.$$ If $x_{_{|M}}=x_{_{|W}}$ then
$$(\chi_{_{W}},\chi_{_M})=\sum_{x\in
\zeta}\frac{|H|}{|\hat{K_g}|}dimM=
\frac{|H||\zeta|}{|\hat{K_g}|}dimM=\frac{|H|}{|K_g|}dimM.$$
This proves part 2 as well, and hence concludes the proof of the
theorem. \qed

\end{document}